\documentclass[eqsecnum,showpacs,preprint,amsmath,aps,nofootinbib]{revtex4}  %

\def\cstok#1{\leavevmode\thinspace\hbox{\vrule height6pt\vtop{\vbox{\hrule\kern-2pt
\hbox{\vphantom{\tt/}\thinspace{\tt#1}\thinspace}}
\kern1pt\hrule}\vrule height6pt}\thinspace}

\begin{document}
\title{A new perspective on the Ermakov-Pinney and scalar wave equations}

\author{Giampiero Esposito ORCID: 0000-0001-5930-8366}
\email[E-mail: ]{gesposit@na.infn.it}
\affiliation{Istituto Nazionale di Fisica Nucleare, Sezione di
Napoli, Complesso Universitario di Monte S. Angelo,
Via Cintia Edificio 6, 80126 Napoli, Italy}

\author{Marica Minucci ORCID: 0000-0002-7095-5115}
\email[E-mail: ]{maricaminucci27@gmail.com}
\affiliation{Dipartimento di Fisica ``Ettore Pancini'', Universit\`a Federico II, \\
Complesso Universitario di Monte S. Angelo, 
Via Cintia Edificio 6, 80126 Napoli, Italy}

\date{\today}

\begin{abstract}
The first part of the paper proves that a subset of the general set of Ermakov-Pinney 
equations can be obtained by differentiation of a first-order non-linear differential
equation. The second part of the paper proves that, similarly, the equation for the
amplitude function for the parametrix of the scalar wave equation can be obtained by
covariant differentiation of a first-order non-linear equation. The construction of
such a first-order non-linear equation relies upon a pair of auxiliary $1$-forms
$(\psi,\rho)$. The $1$-form $\psi$ satisfies the divergenceless condition
${\rm div}(\psi)=0$, whereas the $1$-form $\rho$ fulfills the non-linear equation 
${\rm div}(\rho)+\rho^{2}=0$. The auxiliary $1$-forms $(\psi,\rho)$ are evaluated
explicitly in Kasner space-time, and hence also amplitude and phase function in the
parametrix are obtained. Thus, the novel method developed in this paper can be used
with profit in physical applications.
\end{abstract}


\maketitle

\section{Introduction}

Although the modern theoretical description of gravitational interactions \cite{YCB}
has completely superseded Newtonian gravity, the investigation of ordinary differential
equations provides an invaluable tool in the analysis of chaotic dynamical systems
\cite{Poincare} and in studying the interplay between linear and non-linear
differential equations \cite{CIMM}.

In particular, in our paper we are interested inn
the Ermakov-Pinney \cite{EP1,EP2} non-linear differential equation
\begin{equation}
y''+py=q y^{-3},
\label{(1.1)}
\end{equation}
which has found, along the years, many applications in theoretical physics, including quantum
mechanics \cite{EP3,EP4} and relativistic cosmology \cite{EP5}. The first aim of our work is to
provide yet another perspective on the way of arriving at equations of type (1.1).
For this purpose, section $2$ provides a concise summary of well-established results on
the canonical form of second-order linear differential equations. Section $3$ applies
an ansatz based on amplitude and phase functions, and proves eventually equivalence
between Eqs. of type (1.1) with $p=0$ and our Eq. (3.9), which is a nonlinear equation
with only first derivatives of the desired solution. Section $4$ studies the correspondence 
between sections $2$ and $3$ on the one hand, and the parametrix construction for scalar
wave equation on the other hand. Section $5$ obtains a first-order non-linear equation for 
the amplitude function occurring in such a parametrix. Section $6$ evaluates in 
Kasner space-time the auxiliary $1$-forms that are needed for a successful application 
of our method. Section $7$ solves the first-order equations for amplitude and phase
function in Kasner space-time. Concluding remarks are then made in section $8$.

\section{Canonical form of second-order linear differential equations}

In the theory of ordinary differential equations, it is well-known that every linear
second-order equation
\begin{equation}
\left[{d^{2}\over dx^{2}}+P(x){d \over dx}+Q(x) \right]u(x)=0
\label{(2.1)}
\end{equation}
can be solved by expressing the unknown function $u$ in the form of a product
\begin{equation}
u(x)=\varphi(x) \chi(x),
\label{(2.2)}
\end{equation}
where \cite{WW,E}
\begin{equation}
\varphi(x)={\rm exp} \left(-{1 \over 2} \int P(x)dx \right),
\label{(2.3)}
\end{equation}
while $\chi$ solves the linear equation
\begin{equation}
\left[{d^{2}\over dx^{2}}+J(x) \right]\chi(x)=0,
\label{(2.4)}
\end{equation}
having defined
\begin{equation}
J(x) \equiv Q(x)-{1 \over 4}P^{2}(x)-{1 \over 2}P'(x).
\label{(2.5)}
\end{equation}
All complications arising from the variable nature of coefficient functions $P$ and $Q$ in
Eq. (2.1) are encoded into the potential term $J(x)$ of Eq. (2.4)  
defined in Eq. (2.5). One can therefore
hope to gain insight by the familiar solution of linear second-order equations solved by
$\sin(x),\cos(x)$ or real-valued exponentials. More precisely, a theorem \cite{Valiron} 
guarantees that, if $J(x)$ is continuous on the closed interval $[a,b]$, and if there exist
real constants $\omega,\Omega$ such that 
\begin{equation}
0 < \omega^{2} < J(x) < \Omega^{2},
\label{(2.6)}
\end{equation}
one can compare the zeros of solutions of Eq. (2.4) 
with the zeros of solutions of the equations
\begin{equation}
\eta''+\mu^{2}\eta=0, \;
\mu=\omega \; {\rm or} \; \Omega.
\label{(2.7)}
\end{equation}
Equations (2.7) are solved by periodic functions $\sin(\mu(x-x_{0}))$ which have zeros
at $x_{0}+{k \pi \over \mu}$, $k$ being an integer and $\mu$ being equal to $\omega$
or $\Omega$ as in Eq. (2.7). One can then prove that the difference $\delta$ between two adjacent 
zeros of a solution of Eq. (2.4) lies in between ${\pi \over \Omega}$ and 
${\pi \over \omega}$ \cite{Valiron}.
Equation (2.4) is therefore regarded as the canonical form of every linear second-order
differential equation \cite{WW}.

\section{An ansatz in terms of amplitude and phase functions}

The work in Ref. \cite{Majorana} has shown a long ago that, on looking for
solutions of Eq. (2.4), one can use with profit the ansatz
\begin{equation}
\chi(x)=u(x){\rm exp} \left(i \int \pi(x) u^{-\lambda}(x)dx \right).
\label{(3.1)}
\end{equation}
By doing so, we find
\begin{eqnarray}
\chi''(x)&=& \left[u''(x)+i \pi'(x)u^{1-\lambda}(x)+i(2-\lambda)\pi(x) u'(x)u^{-\lambda}(x)
-\pi^{2}(x) u^{1-2 \lambda}(x) \right] 
\nonumber \\
& \times & {\rm exp} \left(i \int \pi(x) u^{-\lambda}(x)dx \right).
\label{(3.2)}
\end{eqnarray}
If the potential term $J(x)$ vanishes in Eq. (2.4), we therefore find ($A$ and $B$
being integration constants)
\begin{equation}
\chi(x)=A+Bx.
\label{(3.3)}
\end{equation}
On the other hand, by virtue of Eq. (3.3), Eq. (3.2) yields
\begin{equation}
u''(x)+i \pi'(x)u^{1-\lambda}(x)+i(2-\lambda)\pi(x)u'(x)u^{-\lambda}(x)
-\pi^{2}(x)u^{1-2 \lambda}(x)=0.
\label{(3.4)}
\end{equation}
Equation (3.4) suggests setting $\lambda=2$, and hence we find 
\begin{equation}
u^{3}(x)u''(x)-\pi^{2}(x)=-i \pi'(x) u^{2}(x).
\label{(3.5)}
\end{equation}
Thus, if $\pi(x)={\rm constant}=\tau$, we obtain a particular case of the
Ermakov-Pinney equation (1.1) with $p=0$ and $q=\tau^{2}$ therein, i.e.
\begin{equation}
u^{3}(x)u''(x)=\tau^{2}.
\label{(3.6)}
\end{equation}
Furthermore, we can write that
\begin{equation}
\chi(x)=A+Bx=u(x) {\rm exp} \left(i \tau \int 
{dx\over u^{2}(x)} \right),
\label{(3.7)}
\end{equation}
which implies
\begin{equation}
\log(A+Bx)=\log(u(x))+i \tau \int {dx \over u^{2}(x)}.
\label{(3.8)}
\end{equation}
By differentiation, this yields eventually the non-linear equation
\begin{equation}
{u'(x)\over u(x)}+i{\tau \over u^{2}(x)}={B \over (A+Bx)}.
\label{(3.9)}
\end{equation}
In other words, the Ermakov-Pinney equations with $p=0$ in Eq. (1.1) 
are equivalent to the non-linear equation (3.9), {\it provided} that $u$ 
is a function of class $C^{2}$. On the other hand, Eq. (3.9) allows for 
solutions for $u$ which are just of class $C^{1}$.
An useful check of our calculation is obtained by differentiating with respect to $x$
both sides of Eq. (3.9) when $u$ is of class $C^{2}$, 
and then using (3.9) in order to re-express the square
of ${u'(x) \over u(x)}$. One then recovers the Ermakov-Pinney Eq. (3.6), which
therefore originates from Eq. (3.9), and in turn from Eq. (3.8).

\section{Amplitude-phase ansatz for the parametrix of the scalar wave equation}

The work in Ref. \cite{EBD} has studied the parametrix for the scalar wave equation in
curved spacetime. The topic is relevant both for the mathematical theory of hyperbolic equations
on manifolds \cite{E} and for the modern trends in mathematical relativity \cite{M}.
For our purposes, we can limit ourselves to the following outline.

In a pseudo-Riemannian manifold $(M,g)$ endowed with a Levi-Civita connection $\nabla$,
the wave operator
\begin{equation}
\cstok{\ } \equiv \sum_{\mu,\nu=1}^{4}(g^{-1})^{\mu \nu}\nabla_{\nu}\nabla_{\mu}
\label{(4.1)}
\end{equation}
is a variable-coefficient operator, and the homogeneous wave equation $\cstok{\ } \phi=0$, for
given Cauchy data
\begin{equation}
\phi(x,t=0)=\zeta_{0}(x), \; 
{\partial \phi \over \partial t}(x,t=0)=\zeta_{1}(x),
\label{(4.2)}
\end{equation}
can be solved in the form 
\begin{equation}
\phi(x,t)=\sum_{j=0}^{1}E_{j}(t)\zeta_{j}(x),
\label{(4.3)}
\end{equation}
where $E_{j}$ are the Fourier-Maslov integral operators \cite{EBD,T} 
\begin{equation}
E_{j}(t)\zeta_{j}(x)=\sum_{k=1}^{2}(2 \pi)^{-3} \int 
e^{i \varphi_{k}(x,t,\xi)}\alpha_{jk}(x,t,\xi)
{\tilde \zeta}_{j}(\xi)d^{3}\xi
+R_{j}(t)\zeta_{j}(x),
\label{(4.4)}
\end{equation}
the $R_{j}(t)$ being regularizing operators which smooth out the singularities upon
which they act \cite{T}. In the simplest possible terms, the meaning of Eq. (4.4)
is that the integral operators which generalize the Fourier transform to
pseudo-Riemannian manifolds involve again an integrand proportional to
amplitude $\times$ eponential of ($i$ times a phase function), but, unlike
flat space-time, the amplitude function depends explicitly on all
cotangent bundle coordinates, while the phase function is no longer linear
in these variables \cite{T}.

The amplitude and phase functions, denoted by $\alpha$ (of class $C^{2}$) and
$\varphi$ (of class $C^{1}$), respectively, can be obtained by solving the
coupled equations \cite{EBD}
\begin{equation}
\sum_{\gamma,\beta=1}^{4}(g^{-1})^{\gamma \beta}\nabla_{\beta} 
\Bigr(\alpha^{2}\nabla_{\gamma}\varphi \Bigr)=0,
\label{(4.5)}
\end{equation}
\begin{equation}
\sum_{\gamma,\beta=1}^{4} (g^{-1})^{\gamma \beta}
(\nabla_{\beta}\varphi)(\nabla_{\gamma}\varphi)
={\cstok{\ }\alpha \over \alpha}.
\label{(4.6)}
\end{equation}
These equations lead in turn to the following recipe \cite{EBD}. 
First, find a divergenceless covector $\psi_{\gamma}$, i.e.
\begin{equation}
{\rm div}\psi=\sum_{\gamma=1}^{4}\nabla^{\gamma}\psi_{\gamma}
=\sum_{\gamma,\beta=1}^{4} (g^{-1})^{\gamma \beta}
\nabla_{\beta}\psi_{\gamma}=0,
\label{(4.7)}
\end{equation}
then solve the non-linear equation
\begin{equation}
\alpha^{3}\cstok{\ } \alpha=\sum_{\gamma=1}^{4}\psi_{\gamma}\psi^{\gamma}
=\sum_{\gamma,\beta=1}^{4} (g^{-1})^{\gamma \beta}\psi_{\beta}\psi_{\gamma},
\label{(4.8)}
\end{equation}
and eventually obtain the phase from the equation
\begin{equation}
\nabla_{\gamma}\varphi=\alpha^{-2} \psi_{\gamma}.
\label{(4.9)}
\end{equation}

Interestingly, upon defining
\begin{equation}
q \equiv \sum_{\gamma=1}^{4}\psi_{\gamma}\psi^{\gamma},
\label{(4.10)}
\end{equation}
Eq. (4.8) becomes of the type (1.1) with $p=0$. Thus, bearing in mind our finding in
Sec. III, we remark that a simple but non-trivial correspondence exists between a subset
of the general set of Ermakov-Pinney equations and their tensor-calculus counterpart
for the analysis of the scalar wave equation, expressed by the following recipes:
$$
\pi(x)=\tau={\rm constant} \; {\rm in} \; {\rm Eq.} \; (3.5) 
\; \leftrightarrow \; {\rm Eq.} \: (4.7),
$$
$$
{\rm Eq.} \; (3.6) \; \leftrightarrow \; {\rm Eq.} \; (4.8) 
\; {\rm with} \; {\rm constant} \; {\rm value} \; q \;
{\rm of} \; {\rm the} \; {\rm right-hand} \; {\rm side},
$$
$$
u'' \; {\rm in} \; {\rm Eq.} \; (3.6) \; \leftrightarrow \;
\cstok{\ } \alpha \; {\rm in} \; {\rm Eq.} \; (4.8).
$$
Moreover, we know that Eq. (3.6) is solved by a function solving the possibly simpler 
equation (3.9). This implies in turn that Eq. (4.8) for the amplitude $\alpha$ must be
obtainable from the as yet unknown solution ${\cal U}$ of an unknown nonlinear 
equation involving at most first-order derivatives of ${\cal U}$. 
This is the topic of next section.

\section{A first-order non-linear equation for the amplitude in the parametrix}

We are now going to prove that not only does our approach shed new light on the 
Ermakov-Pinney equation as resulting from differentiation of the non-linear equation
(3.9), which is therefore more fundamental (allowing also for solutions which
are only of class $C^{1}$, but not $C^{2}$), but that also the second-order equation 
for the amplitude $\alpha$ in the parametrix can be replaced by a first-order equation.
For this purpose, since $\cstok{\ }\alpha$ should be the counterpart of $u''(x)$, and
the divergenceless condition, Eq. (4.7), the counterpart of $\pi'=0$ in Section $3$, 
we are led to consider the first-order non-linear equation
\begin{equation}
{(\nabla_{\gamma}\alpha)\over \alpha}+i{\psi_{\gamma}\over \alpha^{2}}=\rho_{\gamma},
\label{(5.1)}
\end{equation}
where $\rho_{\gamma}$ are the components of a suitable covector that should generalize the
behaviour of $R(x)={B \over (A+Bx)}$ on the right-hand side of Eq. 
(3.9).\footnote{Note that, strictly speaking, since $\alpha$, $\varphi$ and 
$\psi_{\gamma}$ are real-valued, we are dealing with a complex-valued
vector field $\sum_{\gamma=1}^{4}\rho^{\gamma}{\partial \over \partial x^{\gamma}}$
\cite{OS}, with the associated dual concept of complex-valued $1$-form field.} 
At this stage, inspired by Section $3$, we perform 
covariant differentiation $\nabla^{\gamma}$ of both
sides of Eq. (5.1), finding first the equation
\begin{equation}
-\sum_{\gamma,\beta=1}^{4} (g^{-1})^{\gamma \beta}
{(\nabla_{\beta}\alpha)\over \alpha}
{(\nabla_{\gamma}\alpha)\over \alpha}
+{\cstok{\ }\alpha \over \alpha}
-2i{\psi_{\gamma}\over \alpha^{2}}
{(\nabla^{\gamma}\alpha)\over \alpha}
=\sum_{\gamma,\beta=1}^{4} (g^{-1})^{\gamma \beta}
\nabla_{\beta}\rho_{\gamma},
\label{(5.2)}
\end{equation}
because the divergenceless condition (4.7) holds by assumption. Next, we exploit Eq. (5.1)
by re-expressing all first covariant derivatives of $\alpha$ in Eq. (5.2) in the form
$$
{(\nabla_{\beta}\alpha)\over \alpha}=\rho_{\beta}-i{\psi_{\beta}\over \alpha^{2}},
$$
hence finding
\begin{eqnarray}
\; & \; & -\sum_{\gamma,\beta=1}^{4}(g^{-1})^{\gamma \beta}
\left(\rho_{\beta}-i{\psi_{\beta}\over \alpha^{2}}\right)
\left(\rho_{\gamma}-i{\psi_{\gamma}\over \alpha^{2}}\right)
+{\cstok{\ }\alpha \over \alpha} 
\nonumber \\
&-& 2i \sum_{\gamma=1}^{4}{\psi_{\gamma}\over \alpha^{2}}
\left(\rho^{\gamma}-i{\psi^{\gamma}\over \alpha^{2}}\right)
=\sum_{\gamma,\beta=1}^{4}(g^{-1})^{\gamma \beta}
\nabla_{\beta}\rho_{\gamma}.
\label{(5.3)}
\end{eqnarray}
In this equation, the terms proportional to $\sum_{\gamma=1}^{4}\rho^{\gamma}\psi_{\gamma}$
add up to $0$, and hence we obtain
\begin{equation}
{\cstok{\ }\alpha \over \alpha}
-\sum_{\gamma=1}^{4}{\psi_{\gamma}\psi^{\gamma}\over \alpha^{4}}
=\sum_{\gamma,\beta=1}^{4}(g^{-1})^{\gamma \beta}
\Bigr(\nabla_{\beta}\rho_{\gamma}+\rho_{\beta}\rho_{\gamma} \Bigr).
\label{(5.4)}
\end{equation}
Thus, provided that
\begin{equation}
\sum_{\gamma,\beta=1}^{4}(g^{-1})^{\gamma \beta}
\Bigr(\nabla_{\beta}\rho_{\gamma}+\rho_{\beta}\rho_{\gamma}\Bigr)=0,
\label{(5.5)}
\end{equation}
we obtain eventually the second-order equation (4.8) for the amplitude $\alpha$ in the
parametrix for the scalar wave equation. Remarkably, Eq. (5.5) is precisely the tensorial
generalization of the differential equation obeyed by the right-hand side 
$R(x)={B \over (A+Bx)}$ of Eq. (3.9), because
$$
{d \over dx}R(x)+R^{2}(x)=-B^{2}(A+Bx)^{-2}+B^{2}(A+Bx)^{-2}=0.
$$ 

\section{Evaluation of the auxiliary $1$-forms $\psi$ and $\rho$}

So far, the critical reader might think that our method, despite being elegant and correct, 
does not offer any concrete advantage with respect to the direct investigation of the
coupled equations (4.5) and (4.6), or (4.8) and (4.9). The aim of the present section 
is therefore to prove that the $1$-forms $\psi$ and $\rho$ fulfilling Eqs. (4.7) and 
(5.5) are explicitly computable in a non-trivial case of physical interest.

For this purpose, inspired again by our Ref. \cite{EBD}, we consider Kasner spacetime, 
whose metric in $c=1$ units reads as \cite{YCB}
\begin{equation}
g=-dt \otimes dt+t^{2 p_{1}} dx \otimes dx +t^{2 p_{2}} dy \otimes dy
+t^{2 p_{3}} dz \otimes dz, 
\label{(6.1)}
\end{equation}
where the real numbers $p_{1},p_{2},p_{3}$ satisfy the condition
\begin{equation}
\sum_{k=1}^{3}p_{k}=1,
\label{(6.2)}
\end{equation}
as well as the unit $2$-sphere condition
\begin{equation}
\sum_{k=1}^{3}(p_{k})^{2}=1.
\label{(6.3)}
\end{equation}
Let us assume for simplicity that the only non-vanishing component of the desired 
$1$-form $\psi$ is $\psi_{1}=\psi_{0}(t)$. Hence we find (since the Christoffel 
coefficients $\Gamma_{kk}^{0}=p_{k}t^{2p_{k}-1}$, $\forall k=1,2,3$)
\begin{eqnarray}
{\rm div}(\psi)&=& \sum_{\mu, \nu=1}^{4}(g^{-1})^{\mu \nu}
\nabla_{\nu}\psi_{\mu}
=-\partial_{0}\psi_{0}+\psi_{0}\left(\Gamma_{00}^{0}-\sum_{k=1}^{3}t^{-2p_{k}}
\Gamma_{kk}^{0}\right) 
\nonumber \\
&=& -{d\psi_{0}\over dt}-{\psi_{0}\over t}\sum_{k=1}^{3}p_{k}
=-{d \psi_{0} \over dt}-{\psi_{0}\over t}.
\label{(6.4)}
\end{eqnarray}
The vanishing divergence condition (4.7) is therefore satisfied by
\begin{equation}
\psi_{0}={\kappa \over t}, \;
\kappa={\rm constant}.
\label{(6.5)}
\end{equation}

Similarly, assuming that also the auxiliary $1$-form $\rho$ has only one
non-vanishing component $\rho_{0}(t)$, one finds
\begin{equation}
{\rm div}(\rho)+\rho^{2}=
-{d \rho_{0} \over dt}-{\rho_{0}\over t}+(g^{-1})^{00}(\rho_{0})^{2}
=-{d \rho_{0}\over dt}-{\rho_{0}\over t}-(\rho_{0})^{2}.
\label{(6.6)}
\end{equation}
But we know from section $5$ that $\rho$ should be complex-valued, hence we set
\begin{equation}
\rho_{0}(t)=\beta_{1}(t)+i \beta_{2}(t),
\label{(6.7)}
\end{equation}
$\beta_{1}$ and $\beta_{2}$ being the real and imaginary part of $\rho_{0}$, 
respectively. Thus, by virtue of the identity (6.6), Eq. (5.5) leads to the 
non-linear coupled system
\begin{equation}
{d\beta_{1}\over dt}+{\beta_{1}(t)\over t}
+(\beta_{1}(t))^{2}-(\beta_{2}(t))^{2}=0,
\label{(6.8)}
\end{equation}
\begin{equation}
{d \beta_{2}\over dt}+{\beta_{2}(t)\over t}
+2 \beta_{1}(t) \beta_{2}(t)=0.
\label{(6.9)}
\end{equation}
Equations (6.8) and (6.9) suggest re-expressing them in terms of the 
unknown function
\begin{equation}
B(t) \equiv {\beta_{1}(t) \over \beta_{2}(t)}.
\label{(6.10)}
\end{equation}
This leads to the equivalent system 
\begin{equation}
{dB \over dt}\beta_{2}(t)+B(t) \left({d \beta_{2}\over dt}+
{\beta_{2}(t)\over t}\right)+(B^{2}(t)-1)(\beta_{2}(t))^{2}=0,
\label{(6.11)}
\end{equation}
\begin{equation}
{d \beta_{2}\over dt}+{\beta_{2}(t)\over t}
+2B(t)(\beta_{2}(t))^{2}=0.
\label{(6.12)}
\end{equation}
By insertion of Eq. (6.12) into Eq. (6.11), we find
\begin{equation}
\beta_{2}(t)={B'(t)\over (1+B^{2}(t))},
\label{(6.13)}
\end{equation}
\begin{equation}
\beta_{1}(t)=B(t)\beta_{2}(t)={B(t)B'(t)\over (1+B^{2}(t))},
\label{(6.14)}
\end{equation}
and hence Eq. (6.12) yields for $B(t)$ the equation
\begin{equation}
{B''(t)\over (1+B^{2}(t))}-2{B(t)(B'(t))^{2}\over (1+B^{2}(t))^{2}}
+{1 \over t}{B'(t)\over (1+B^{2}(t))}
+2{B(t)(B'(t))^{2}\over (1+B^{2}(t))^{2}}=0,
\label{(6.15)}
\end{equation}
which is equivalent to the linear differential
equation\footnote{Note also that $(1+B^{2}(t))$ in Eq. (6.15) can never vanish, bearing 
in mind the real nature of $B(t)$ from the definition (6.10).}
\begin{equation}
B''(t)+{1 \over t}B'(t)=0.
\label{(6.16)}
\end{equation}
Equation (6.16) implies that $B'(t)$ is proportional to ${1 \over t}$, and hence,
upon introducing the real parameter $\sigma$, one can write that 
($\kappa$ being the same parameter used in (6.5))
\begin{equation}
{dB \over dt}={\sigma \over \kappa}{1 \over t},
\label{(6.17)}
\end{equation}
and hence
\begin{equation}
B(t)=B(T)+{\sigma \over \kappa}\log \left({t \over T} \right).
\label{(6.18)}
\end{equation}
Hereafter we set $B(T)=0$ for simplicity. By virtue of (6.13), (6.14) and (6.18)
we obtain eventually, upon defining\footnote{The work in Ref. \cite{BE19} 
arrives instead at Eq. (6.19)
by solving directly for the amplitude $\alpha$, without making any use 
of the auxiliary $1$-form $\rho$ and of our Eq. (5.1).}
\begin{equation}
D_{\sigma}(t) \equiv {\kappa^{2}\over \sigma}
+\sigma \log^{2} \left({t \over T}\right),
\label{(6.19)}
\end{equation}
the exact formulae
\begin{equation}
\beta_{1}(t)={\sigma \over t}
{\log \left({t \over T}\right) \over D_{\sigma}(t)},
\label{(6.20)}
\end{equation}
\begin{equation}
\beta_{2}(t)={\kappa \over t}{1 \over D_{\sigma}(t)}.
\label{(6.21)}
\end{equation}
Remarkably, the exact solution of the non-linear equations (6.8) and (6.9) has been
obtained from the general solution of the linear equation (6.16).

\section{Amplitude and phase functions in Kasner space-time}

In light of (5.1), we can now evaluate the amplitude function
$\alpha$ from the equation 
\begin{equation}
\nabla_{\gamma}(\log(\alpha))+i{\psi_{\gamma}\over \alpha^{2}}=\rho_{\gamma},
\label{(7.1)}
\end{equation}
and eventually the phase function $\varphi$ from Eq. (4.9), which reads
in our case
\begin{equation}
{d \varphi \over dt}={\kappa \over t \alpha^{2}}.
\label{(7.2)}
\end{equation}
From Eq. (7.1) we obtain, in Kasner space-time, the ordinary differential equation
\begin{equation}
{d \over dt}\log (\alpha(t))+i{\kappa \over t \alpha^{2}(t)}
=\beta_{1}(t)+i \beta_{2}(t),
\label{(7.3)}
\end{equation}
i.e., upon separating real and imaginary part, the pair of equations
\begin{equation}
{d \over dt}\log(\alpha(t))=\beta_{1}(t),
\label{(7.4)}
\end{equation}
\begin{equation}
{\kappa \over t \alpha^{2}(t)}=\beta_{2}(t)=\varphi(t).
\label{(7.5)}
\end{equation}
Hence we find in Kasner space-time the amplitude function
\begin{equation}
\alpha(t)=\alpha(T){\rm exp} \int_{T}^{t}
{\sigma \over \tau}{\log \left({\tau \over T}\right) \over
D_{\sigma}(\tau)}d\tau
=\sqrt{{\kappa^{2}\over \sigma}+\sigma \log^{2}
\left({t \over T}\right)},
\label{(7.6)}
\end{equation}
for which $\alpha(T)=\sqrt{{\kappa^{2}\over \sigma}}$,
as well as the phase function
\begin{equation}
\varphi(t)= \varphi(T)+\kappa 
\int_{T}^{t}{d \tau \over \tau D_{\sigma}(\tau)}
= \varphi(T)+{\rm arctan} \left({\sigma \over \kappa}
\log \left({t \over T}\right)\right),
\label{(7.7)}
\end{equation}
which holds for all positive values of the real ratio ${T \over (t-T)}$.
It should be stressed that, in a generic space-time without any symmetry,
the amplitude and phase, if computable, will depend on all cotangent bundle
local coordinates \cite{T} (see further comments in Section $8$).

\section{Concluding remarks}

In our paper, starting from well known properties in the theory of linear differential
equations, we have first proved that the Ermakov-Pinney equations with $p=0$ in Eq. 
(1.1) result from differentiation of the more fundamental equation (3.9), provided
that the function $u$ solving (3.9) is taken to be at least of class $C^{2}$.

By comparison with the construction of amplitude and phase in the scalar parametrix, we 
have then proved that finding the amplitude $\alpha$ for which Eq. (4.8) holds
with a divergenceless covector $\psi_{\gamma}$, is equivalent to finding also a
covector $\rho_{\gamma}$ for which Eq. (5.5) holds. One can then obtain the amplitude
$\alpha$ from the first-order non-linear equation (7.1). 
Our successful calculations of sections $6$ and $7$, where we have evaluated the 
auxiliary $1$-forms with components
\begin{equation}
\psi_{\mu}(t)=\left({\kappa \over t},0,0,0 \right),
\label{(8.1)}
\end{equation}
\begin{equation}
\rho_{\mu}(t)=\left({1 \over t D_{\sigma}(t)}
\left(\sigma \log \left({t \over T}\right)+i \kappa \right),0,0,0 \right),
\label{(8.2)}
\end{equation}
prove that our original method leads to a powerful tool for 
studying the scalar wave equation with the associated 
parametrix. This will be of concrete interest in applied mathematics and in the
theoretical physics of fundamental interactions.

Note also that, in principle, there might exist solutions of Eq. (7.1) which are of
class $C^{1}$ but not $C^{2}$. Thus, the consideration of Eq. (5.1) is closer to the
modern emphasis on finding new solutions of partial differential equations under
weaker differentiability properties. Of course, the corresponding physical interpretation
is a relevant open problem. 

\acknowledgments
G. E. is grateful to Dipartimento di Fisica ``Ettore Pancini'' for hospitality and support,
and to D. Bini for a careful reading of a draft version, and for encouragement.


\begin{thebibliography}{}

\bibitem{YCB}
Y. Choquet-Bruhat, {\it General Relativity and the Einstein Equations}
(Clarendon Press, Oxford, 2009).

\bibitem{Poincare}
H. Poincar\'e, {\it Les M\'ethodes Nouvelles de la M\'ecanique C\'eleste}. Tomes 1, 2, 3
(Gauthier-Villars, Paris, 1892; ibid. 1893; ibid. 1899).

\bibitem{CIMM}
J. Carinena, A. Ibort, G. Marmo and G. Morandi, 
{\it Geometry from Dynamics, Classical and Quantum} (Springer, Berlin, 2015).

\bibitem{EP1}
V. P. Ermakov, Univ. Izv. Kiev Ser. III {\bf 9}, 1 (1880).

\bibitem{EP2}
E. Pinney, The nonlinear differential equation $y''+p(x)y+c y^{-3}=0$, 
Proc. Amer. Math. Soc. {\bf 1}, 681 (1950).

\bibitem{EP3}
H. R. Lewis and W. B. Riesenfeld, An exact quantum theory of the time-dependent
harmonic oscillator and of a charged particle in a time-dependent electromagnetic 
field, J. Math. Phys. {\bf 10}, 1458 (1969).

\bibitem{EP4}
P. G. Kevrekidis and Y. Drossimos, Nonlinearity from linearity: the Ermakov-Pinney equation
revisited, Math. Comp. Sim. {\bf 74}, 196 (2007).

\bibitem{EP5}
R. M. Hawkins and J. E. Lidsey, Ermakov-Pinney equation in scalar-field cosmologies, 
Phys. Rev. D {\bf 66}, 023523 (2002).

\bibitem{WW}
E. T. Whittaker and G. N. Watson, {\it Modern Analysis}
(Cambridge University Press, Cambridge, 1927).

\bibitem{E}
G. Esposito, {\it From Ordinary to Partial Differential Equations}, Unitext {\bf 106}
(Springer, Berlin, 2017).

\bibitem{Valiron}
G. Valiron, {\it \'Equations Fonctionnelles. Applications} (Masson, Paris, 1945).

\bibitem{Majorana}
E. Majorana, {\it Notes on Theoretical Physics}, edited by
S. Esposito and E. Recami (Springer, Berlin, 2003).

\bibitem{EBD}
G. Esposito, E. Battista and E. Di Grezia, Bicharacteristics and Fourier integral
operators in Kasner spacetime, Int. J. Geom. Methods Mod. Phys. {\bf 12}, 1550060 (2015).

\bibitem{M}
M. Minucci, {\it Hyperbolic Equations and General Relativity}
(Nova Science, New York, 2019).

\bibitem{T}
F. Treves, {\it Introduction to Pseudodifferential and Fourier Integral Operators.
Fourier Integral Operators}, Vol. 2 (Plenum Press, New York, 1980). 

\bibitem{OS}
O. Stormark, On the theorem of Frobenius for complex vector fields, Ann. Sc. Norm.
Sup. Serie $4$, {\bf 9}, 57 (1982).

\bibitem{BE19}
D. Bini and G. Esposito, Ermakov-Pinney-like non-linear equation in curved space-time,
work in progress (2019).

\end{thebibliography}
\end{document}